\documentclass{article}\begin{document}
\centerline{\bf Properties of Lerch Sums and Ramanujan's Mock Theta Functions}

\[
\]

\centerline{Nikos Bagis}

\centerline{Stenimahou 5 Edessa}
\centerline{Pella 58200, Greece}
\centerline{nikosbagis@hotmail.gr}

\[
\]

\centerline{\bf Abstract}

\begin{quote}
In this article we study properties of Ramanujan's mock theta functions that can be expressed by Lerch sums. We mainly show that each Lerch sum is actually the integral of a Jacobian theta function, (here we show that for $\vartheta_3(t,q)$ and $\vartheta_4(t,q)$) multiplied by the $\sec-$function ($\textrm{secant}$). We also prove some modular relations and evaluate the Fourier coefficients of a class of Lerch sums.          
\end{quote}


\section{Introduction (Lerch sums and Ramanujan's mock theta functions)}

We first introduce the $q-$Pochhammer symbol
\begin{equation}
(a;q)_n:=\prod^{n}_{j=1}\left(1-aq^{j-1}\right)\textrm{, }n=1,2,\ldots,
\end{equation}
where $(a;q)_0:=1$ and start with some remarks regarding the mock theta function  
\begin{equation}
f(q)=\sum^{\infty}_{n=0}\frac{q^{n^2}}{(-q;q)^2_n}.
\end{equation}
In parallel we will expose the main idea of Zwegers work [6]. 

The function $f(q)$ is an order 3 mock theta function. If we define the Ramanujan $\eta-$eta function by
\begin{equation}
\eta(q):=\prod^{\infty}_{n=1}\left(1-q^n\right)\textrm{, }|q|<1,
\end{equation}
then we have a first known result due to Watson.\\ 
\\
\textbf{Theorem 1.}
\begin{equation}
f(q)=\frac{2}{\eta(q)}\sum^{\infty}_{n=-\infty}\frac{(-1)^nq^{3n^2/2+n/2}}{1+q^n}\textrm{, }|q|<1.
\end{equation}
\\

For a list of more mock theta functions regarding those given by Ramanujan see:\\
\\ 
\centerline{https://en.wikipedia.org/wiki/Mock modular form}\\
\\
An infinite class of mock theta functions can expressed in the form of Lerch sums. Actually, Zwegers in his doctoral thesis [6] constructed a systematic way to recover them and produce new mock theta functions from Lerch sums. 

According to Zwegers, if we multiply $f(q)$ with a rational power of $q$, here $-1/24$ to get $h(z)=q^{-1/24}f(q)$, $q=e(z):=e^{2\pi i z}$ and then add the non holomorphic function (unary theta series):
\begin{equation}
R(z)=\sum_{\scriptsize 
\begin{array}{cc}
	n\in \textbf{Z}\\
	n\equiv 1(6)
\end{array}
\normalsize}\textrm{sign}(n)\beta\left(n^2y/6\right)q^{-n^2/24}\textrm{, }y=Im(z), 
\end{equation}
where 
\begin{equation}
\beta(x)=\int^{\infty}_{x}t^{-1/2}e^{-\pi t}dt\textrm{, }x\geq 0,
\end{equation}
then, the result $h(z)+R(z)$ becomes a modular form of weight $1/2$ in the modular group $\Gamma(2)$.

More precisely and generally Zwegers showed, that for all  Lerch sums (normalized Lerch sums): 
\begin{equation}
\mu(u,v;\tau)=\frac{a^{1/2}}{\theta(v;\tau)}\sum^{\infty}_{n=-\infty}\frac{(-b)^nq^{n(n+1)/2}}{1-aq^n},
\end{equation}
where $q=e(\tau)$, $a=e(u)$, $b=e(v)$, $Im(\tau)>0$ and 
\begin{equation}
\theta(v;\tau)=\sum_{\scriptsize
\begin{array}{cc}
	\nu\in\textbf{Z}+\frac{1}{2}
\end{array}
\normalsize}(-1)^{\nu-1/2}q^{\nu^2/2}b^{\nu},
\end{equation}
there exists a non holomorphic function
\begin{equation}
R(z;\tau)=\sum_{\scriptsize 
\begin{array}{cc}
	\nu\in\textbf{Z}+\frac{1}{2}
\end{array}
\normalsize}(-1)^{\nu-1/2}\left[\textrm{sign}(\nu)-E\left(\left(\nu+\frac{Im(z)}{y}\right)\sqrt{2y}\right)\right]e(-\nu z)q^{-\nu^2/2},
\end{equation}
where $y=Im(\tau)$ and
\begin{equation}
E(z)=\textrm{erf}(\sqrt{\pi}z)=2\int^{z}_{0}e^{-\pi u^2}du,
\end{equation}
such that the function
\begin{equation}
M(u,v;\tau)=\mu(u,v;\tau)-\frac{1}{2}R(u-v;\tau),
\end{equation}
becomes meromorphic Jacobi form and the normalized Lerch sum (7) is a mock theta function. In our case (function (2)), Ramanujan's $f(q)$ function is a mock theta function.

Actualy all mock theta functions of Ramanujan are the holomorphic part of a weight $1/2$ harmonic weak Maass forms. This fact have led to interesting connections of mock theta functions with number theory (see [10],[11]). Also Ramanujan's partial theta functions (mock theta functions) have relation with the construction of higher dimensional multisums, (see [7] and [12]). For applications of mock theta functions in modern physics and quantum black holes a representative paper is that of A. Dadholkar, S. Murthy and D. Zagier [13].

Also recently, Hickerson and Mortenson have expressed all mock theta functions in terms of Appell-Lerch sums (see [8],[9]).\\

In this paper we examine of how we can treat with general Lerch sums like ($q=e(z)$, $Im(z)>0$):
$$
\sum^{\infty}_{n=-\infty}\frac{q^{an^2+bn}}{\cosh(2\pi w(An+B))}.
$$
This sum looks like a theta function only in numerator. However as we will show its relation with a classical Jacobian theta function $\vartheta_3(z,q)$ is (when $a,A,Im(z)>0$, $Im(w)\neq0$):     
$$
\sum^{\infty}_{n=-\infty}\frac{q^{an^2+bn}}{\cosh(2\pi i w(An+B))}=
$$
$$
=\frac{e^{-2iBb\pi z/A}}{2A\pi i w}\int^{+\infty}_{-\infty}\vartheta_3\left(\gamma,e(az)\right)e^{-2iB\gamma/A}\sec\left(\frac{\gamma+b\pi z}{2Aw}\right)d\gamma.
$$
This integral expansion of the Lerch sums have very intersting interpretations theoretical and practical. It enable us to evaluate the Lerch sums in terms of Jacobi theta functions. The Jacobi theta functions have known properties, hence we can extract results for Ramanujan's mock theta functions using the transformation properties of Jacobi theta functions. Also the nature of the sums give us litle information about the function to which they converge. Now if we wish to have a ''nice'' transformation which overclose all of them (Lerch sums and theta functions) we can set a Fourier exponential inside the integral and generalize the Lerch sums. The applications for the Lerch sums can easily obtained if we set the exponetntial to be 1. Proceeding in this way we find indeed interesting transformation formulas for the Lerch sums from that of Jacobi theta functions. In the rest of the paper we give evaluations for the Fourier coefficients of the mock theta functions using divisor sums and give examples. Lastly, we give some miscellaneous results for the representations of the $f(q)$ function and a theorem for the transformation of theta functions, presented and proved in a previous work, but not having much generality.

\section{Some remarks about $f(q)$}

We now return to (2) to examine its convergence regions and some enalactic expresions. 

The next theorem concerns about the range of convergence of $f(q)$ when $q\in\textbf{C}$ and as we will see $f(q)$ (given by (2)) converges in the whole complex plane with exception the unit circle $|q|=1$
For to prove these we first need a lemma.\\
\\ 
\textbf{Lemma 1.}\\
i) If $|q_1|>1$ and $a\neq q_1^{-n-k+1}$, $n,k\in\textbf{N}$, then it holds
\begin{equation}
(a;q_1)_n=\frac{\left(\frac{1}{a};\frac{1}{q_1}\right)_{\infty}}{\left(\frac{1}{aq_1^{n}};\frac{1}{q_1}\right)_{\infty}}\frac{1}{(-1/a)^nq_1^{-n(n-1)/2}}.
\end{equation}
ii) If $|q|<1$ and $a\neq q^{-n-k+1}$, $n,k\in\textbf{N}$, then it holds
\begin{equation}
(a;q)_n=\frac{(a;q)_{\infty}}{(aq^n;q)_{\infty}}.
\end{equation}
\\
\textbf{Proof.}\\
We first prove (13) and then (12). Under the restrictions for $a$ in (13), we have  
$$
\frac{(a;q)_{\infty}}{(aq^n;q)_{\infty}}=\frac{\prod^{\infty}_{k=1}\left(1-aq^{k-1}\right)}{\prod^{\infty}_{k=1}\left(1-aq^{n+k-1}\right)}=\prod^{n}_{k=1}\left(1-aq^{k-1}\right)=(a;q)_n.
$$
For (12) we write $q_1=1/q$, $|q|<1$ and using (13) we have
$$
\frac{\left(\frac{1}{a};\frac{1}{q_1}\right)_{\infty}}{\left(\frac{1}{aq_1^n};\frac{1}{q_1}\right)_{\infty}}=\left(\frac{1}{a},\frac{1}{q_1}\right)_n=\prod^{n}_{k=1}\left(1-\frac{1}{a}\frac{1}{q_1^{k-1}}\right)=\prod^{n}_{k=1}\frac{1}{aq_1^{k-1}}\left(aq_1^{k-1}-1\right)=
$$
$$
=(a;q_1)_n\frac{(-1/a)^n}{q_1^{n(n-1)/2}}.
$$
$qed$\\
\\
\textbf{Theorem 2.}\\
The function $f(q)$ is well defined for all $q\in\textbf{C}-D$, where $D=\{z\in\textbf{C}:|z|=1\}$. Moreover if $|q|<1$, then
\begin{equation}
f(q)=\sum^{\infty}_{n=0}\frac{q^{n^2}}{(-q;q)_n^2},
\end{equation}
\begin{equation}
f(1/q)=\sum^{\infty}_{n=0}\frac{q^{n}}{(-q;q)_n^2}
\end{equation}
and
\begin{equation}
f(q)=\chi(q)^{-2}\sum^{\infty}_{n=0}\left(-q^{n+1};q\right)_{\infty}^2q^{n^2},
\end{equation}
\begin{equation}
f(1/q)=\chi(q)^{-2}\sum^{\infty}_{n=0}\left(-q^{n+1};q\right)_{\infty}^2q^{n}
\end{equation}
and
$$
f(q)=\chi(q)^{-2}\sum^{\infty}_{n=0}q^{n^2}\exp\left[-2\sum^{\infty}_{s=1}\frac{q^{s}}{s}\sum_{d|s,d\geq n+1}(-1)^{s/d}d\right]
$$
and
\begin{equation}
f(q)\chi(q)^{2}=\sum^{\infty}_{n=0}q^{n^2}\exp\left[-2\sum^{\infty}_{s=1}q^{s}\sum_{\scriptsize
\begin{array}{cc}
0<d|s\\	
d\leq s/(n+1)
\end{array}
\normalsize}\frac{(-1)^{d}}{d}\right].
\end{equation}
Also
\begin{equation}
f(1/q)\chi(q)^{2}=\sum^{\infty}_{n=0}q^n\exp\left[-2\sum^{\infty}_{s=1}q^{s}\sum_{\scriptsize
\begin{array}{cc}
0<d|s\\	
d\leq s/(n+1)
\end{array}
\normalsize}\frac{(-1)^{d}}{d}\right],
\end{equation}
where $\chi(q)=(-q,q)_{\infty}$.\\
\\
\textbf{Proof.}\\
Set $q_1=1/q$. Then for to prove (15) we use 
\begin{equation}
\left(-q_1;q_1\right)_n=\prod^{n}_{j=1}\left(1+q_1^j\right)=\prod^{n}_{j=1}q_1^j\left(1+q_1^{-j}\right)=q_1^{n(n+1)/2}\left(-q;q\right)_n.
\end{equation}
For proving (17) we use (12) of Lemma 1. We have
\begin{equation}
(-q_1;q_1)_{n}=\frac{(-q;q)_{\infty}}{(-q^{n+1};q)_{\infty}}q^{-n(n+1)/2}.
\end{equation}
For to prove (16) we use (20) and (21) in (14).  
For to prove (18) and (19) we write (for (19)) 
$$
f(q_1)=\sum^{\infty}_{n=0}\frac{\left(-q^{n+1};q\right)_{\infty}^2}{\left(-q;q\right)_{\infty}^2}q^n
=\chi(q)^{-2}\sum^{\infty}_{n=0}q^n\prod^{\infty}_{m=1}\left(1+q^nq^m\right)^2=
$$
$$
=\chi(q)^{-2}\sum^{\infty}_{n=0}q^n\exp\left[2\log\left(\prod^{\infty}_{m=1}\left(1+q^nq^m\right)\right)\right]=
$$
$$
=\chi(q)^{-2}\sum^{\infty}_{n=0}q^n\exp\left[2\sum^{\infty}_{m=1}\log(1+q^nq^m)\right]=
$$
$$
=\chi(q)^{-2}\sum^{\infty}_{n=0}q^n\exp\left[-2\sum^{\infty}_{m=1}\sum^{\infty}_{l=1}q^{nl}q^{ml}\frac{(-1)^l}{l}\right]=
$$
$$
=\chi(q)^{-2}\sum^{\infty}_{n=0}q^n\exp\left[-2\sum^{\infty}_{m,l=1}q^{(n+m)l}\frac{(-1)^l}{l}\right]=
$$
$$
=\chi(q)^{-2}\sum^{\infty}_{n=0}q^n\exp\left[-2\sum^{\infty}_{s=1}q^{s}\sum_{d|s,d\geq n+1}\frac{(-1)^{s/d}}{s/d}\right].
$$
Hence finally
$$
f(1/q)=\chi(q)^{-2}\sum^{\infty}_{n=0}q^n\exp\left[-2\sum^{\infty}_{s=1}\frac{q^{s}}{s}\sum_{d|s,d\geq n+1}(-1)^{s/d}d\right].
$$
In the same way we get (18). $qed$

\section{Representations and Fourier coefficients of some general class of Lerch sums}

Theorem 3 below will help us link Lerch sums with the  classical Jacobi theta functions. As we will see it turns out that some Lerch sums are just the image of an integral transform of a Jacobian theta function.\\
\\
\textbf{Theorem 3.}\\
Let $q$ be complex number such $|q|<1$ and $g$ be function with Fourier series expansion
\begin{equation}
g(\phi,q)=\sum^{\infty}_{n=1}a_n(q)\cos(2n\phi),
\end{equation}
then
\begin{equation}
\int^{\pi}_{0}g(\phi,q)\log\left(\frac{\vartheta_4(\phi,q)}{\vartheta_4(0,q)}\right)d\phi=-\pi\sum^{\infty}_{n=1}\frac{a_n(q)q^n}{n(1-q^{2n})},
\end{equation}
where
\begin{equation}
\vartheta_4(z,q)=\sum^{\infty}_{n=-\infty}(-1)^n q^{n^2}e^{2 i n z}
\end{equation}
\\
\textbf{Proof.}\\
See [2].\\
\\
\textbf{Proposition 1.}\\
If $z,q\in\textbf{C}$, with $|q|<1$ we define
\begin{equation}
\psi(z,q):=\sum^{\infty}_{n=-\infty}(-1)^n q^{3n^2/2+n/2}\sin(2zn).
\end{equation}
Then
$$
1+\frac{2}{\pi}\int^{\pi}_{0}\partial_{t}\psi(t,q)\log\left(\frac{\vartheta_4(t,q)}{\vartheta_4(0,q)}\right)dt=\eta(q)f(q)=
$$
\begin{equation}
=2\sum^{\infty}_{n=-\infty}\frac{(-1)^nq^{3n^2/2+n/2}}{1+q^n}
\end{equation}
\\
\textbf{Proof.}\\
Proposition 1 is direct application of Theorem 3.\\
\\

Continuing inspired from Proposition 1 we give the next  generalized theorem:\\ 
\\
\textbf{Theorem 4.}\\
If $|q|<1$ and $a>0$ and if
\begin{equation}
\psi_1(t,q)=\sum^{\infty}_{n=1}(-1)^n n q^{an^2+bn}\cos(2nt),
\end{equation}
then
\begin{equation}
\int^{\pi}_{0}\psi_1(t,q)\log\left(\frac{\vartheta_4(t,q)}{\vartheta_4(0,q)}\right)dt=-\pi\sum^{\infty}_{n=1}\frac{(-1)^nq^{an^2+(b+1)n}}{1-q^{2n}}
\end{equation}
Also if
\begin{equation}
\psi_2(t,q)=-\sum^{\infty}_{n=1}(-1)^n n q^{an^2-bn}\cos(2nt),
\end{equation}
then
\begin{equation}
\int^{\pi}_{0}\psi_2(t,q)\log\left(\frac{\vartheta_4(t,q)}{\vartheta_4(0,q)}\right)dt=\pi\sum^{\infty}_{n=1}\frac{(-1)^nq^{an^2+(-b+1)n}}{1-q^{2n}}
\end{equation}
Hence if
\begin{equation}
\psi(t,q)=\psi_1(t,q)+\psi_2(t,q)=\sum^{\infty}_{n=-\infty}(-1)^nq^{an^2+bn}n\cos(2nt),
\end{equation}
we get
\begin{equation}
I=\int^{\pi}_{0}\psi(t,q)\log\left(\frac{\vartheta_4(t,q)}{\vartheta_4(0,q)}\right)dt
=\pi \sum_{n\in\textbf{\scriptsize Z\normalsize }-\{0\}}\frac{(-1)^n q^{an^2+bn}}{q^{-n}-q^n}
\end{equation}
and if $q=e(z)$, $Im(z)>0$, then
\begin{equation}
\frac{I}{\pi}=\frac{1}{2}\sum_{n\in\textbf{\scriptsize Z}-\{0\}}\frac{(-1)^n\exp\left(2\pi i(an^2+bn)z\right)}{\sinh(2\pi i n z)}.
\end{equation}
\\

The next theorem is representation (reformulation) and evaluation of certain cases of Lerch sums using theta functions.\\
\\
\textbf{Theorem 5.}\\
For $a,b,c\in\textbf{R}$, $a>0$, $c\neq0$, $q=e(z)$, $Im(z)>0$, we define the function 
\begin{equation}
f_s(a,b;c;z):=\sum_{
n\in\textbf{\scriptsize Z\normalsize}-\{0\}}\frac{(-1)^ne\left((an^2+bn)z\right)}{\sinh(2\pi i n z c)}.
\end{equation}
Then:\\
1) If $a$ positive real, $b\in \textbf{N}$, $c=1$, the function $f_s(a,b;c;z)$ is a finite sum of ordinary theta functions i.e
\begin{equation}
f_s(a,b;c;z)=\sum_{j}'\left(2\sum^{\infty}_{n=1}(-1)^n e(an^2z)\cos(b_j\pi n z)\right),
\end{equation}
where the $b_j$ are given from relations (38),(39) below.\\
2) When $a,b,c$ are real and $a>0$, $c\neq 0$, then
$$
f_s(a,b;c;z)=\sum^{\infty}_{l=0}\sum^{\infty}_{n=1}(-1)^n q^{a n^2+c(2l+1)n-bn}-\sum^{\infty}_{l=0}\sum^{\infty}_{n=1}(-1)^n q^{an^2+c(2l+1)n+bn}=
$$
\begin{equation}
=\sum_{l,n\in\textbf{\scriptsize Z\normalsize}}(-1)^n\epsilon(n,l)q^{an^2+(c-b)n+2cln}+\sum^{\infty}_{n=1}(-1)^nq^{an^2+(c-b)n},
\end{equation}
where  
\begin{equation}
\epsilon(n,l):=
\textrm{sign}(n)\textrm{sign}(l)\frac{\textrm{sign}(n)+\textrm{sign}(l)}{2}.
\end{equation}
\\
\textbf{Proof.}\\
1) If $b\in\textbf{N}$ then the next expansion holds
\begin{equation}
\frac{1}{2}\sum_{n\in\textbf{\scriptsize Z}-\{0\}}\frac{(-1)^ne\left((an^2+bn)z\right)}{\sinh(2\pi i n z)}=\sum^{\infty}_{n=1}(-1)^ne(an^2z)U_{b-1}(\cos(2\pi nz)),
\end{equation}
where $U_{n}(x)$ is the $n-$th degree Chebyshev polynomial. Hence
\begin{equation}
U_{b-1}(\cos(2\pi n z))=\left\{
\begin{array}{cc}
1+2\sum^{\left[\frac{b-1}{2}\right]}_{j=1}\cos(4j\pi n z)\textrm{, when }b-\textrm{odd}\\
2\sum^{\left[\frac{b-1}{2}\right]}_{j=0}\cos((4j+2)\pi nz)\textrm{, when }b-\textrm{even}.	
\end{array}\right\}
\end{equation}
and the first result follows.\\
2) The first equality of (36) is obtained easily from  
\begin{equation}
\frac{e(bnz)-e(-bnz)}{2\sinh(2\pi i n z c)}=-(e(bnz)-e(-bnz))\sum^{\infty}_{l=0}e(cn(2l+1)z).
\end{equation} 
The second equality of (36) is obtained using the symbols $\epsilon(n,l)$. $qed$\\
\\
\textbf{Theorem 6.}\\
If $a,b,c\in\textbf{Z}$ and $a>0$, $c\neq 0$ and $q=e(z)$, $Im(z)>0$, then the function
\begin{equation}
f_s(a,b;c;z)=\sum_{n\in\textbf{\scriptsize Z}-\{0\}}\frac{(-1)^ne\left((an^2+bn)z\right)}{\sinh(2\pi i n c z)},
\end{equation}
have the following property:\\
If $n$ is positive integer we define 
\begin{equation}
C_1(a,b,c;n):=2\sum_{\scriptsize
\begin{array}{cc}
\textrm{abs}(d)|n\\
d\neq0	
\end{array}
\normalsize}(-1)^dX_{\textbf{\scriptsize Z}}\left(\frac{b+ad-n/d}{c}\right)\epsilon\left(d,-\frac{b+ad-n/d}{c}\right),
\end{equation}
where $X_{\textbf{\scriptsize{Z}}}(.)$ is the characteristic on integers. Then if we set
$$
C_s(a,b;c;n)=C_1(a,c-b;2c;n),
$$
we have
\begin{equation}
f_s\left(a,b;c;z\right)=\sum^{\infty}_{n=1}C_s(a,b;c;n)q^n+2\sum^{\infty}_{n=1}(-1)^nq^{a n^2+(c-b) n}\textrm{, }|q|<1.
\end{equation}
\\
\textbf{Proof.}\\
Assume the form
\begin{equation}
Q=n\left(an+b+cm\right).
\end{equation}
Using the transformation
\begin{equation}
a=k_1A+l_1B\textrm{, }c=k_2A+l_2B\textrm{, }A,B\in\textbf{Z}-\{0\}\textrm{, }(A,B)=1,
\end{equation}
the form $Q$ becomes
$$
Q=n\left(n(k_1A+l_1B)+b+m(k_2A+l_2B)\right)=
$$
\begin{equation}
=n\left(A(nk_1+mk_2)+b+B(nl_1+ml_2)\right).
\end{equation} 
If 
$$
nk_1+mk_2=n'\textrm{, }nl_1+ml_2=m'
$$
and if  
\begin{equation}
k_1l_2-l_1k_2=1\textrm{, }k_1,k_2,l_1,l_2,\in \textbf{Z}
\end{equation} 
we get
$$
Q=n\left(An'+b+Bm'\right)
$$
and $n=n'l_2-m'k_2$. By this way $Q$ becomes 
\begin{equation}
Q=\left(n'l_2-m'k_2\right)\left(An'+Bm'+b\right).
\end{equation}
Hence if we set in $Q$
\begin{equation}
n'_1=n'l_2-m'k_2\textrm{, }n'_2=An'+Bm'+b,
\end{equation}
then
\begin{equation}
Q=n'_1n'_2.
\end{equation}
This case corresponds to $a,b,c\in\textbf{Z}$, $a>0$. Going backwards we get 
\begin{equation}
n=n'_1\textrm{, }m=-\frac{b+an'_1-n'_2}{c},
\end{equation}
and $Q$ becomes 
\begin{equation}
Q_1=n'_1n'_2.
\end{equation}
Hence we obtain the next identity
$$
g\left(q\right)=\sum^{\infty}_{n,m=-\infty}(-1)^n\epsilon(n,m)q^{n(an+b+cm)}=
$$
$$
=\sum_{\scriptsize
\begin{array}{cc}
n_1',n_2'\in\textbf{Z}\normalsize
\end{array}
}(-1)^{n_1'}X_{\textbf{\scriptsize{Z}}}\left(\frac{b+an'_1-n'_2}{c}\right)\epsilon\left(n'_1,-\frac{b+an'_1-n'_2}{c}\right)q^{n'_1n'_2}=
$$
\begin{equation}
=\sum^{\infty}_{n=1}\left(\sum_{\scriptsize
\begin{array}{cc}
\textrm{abs}(d)|n\\
d\neq0	
\end{array}
\normalsize}(-1)^dX_{\textbf{\scriptsize{Z}}}\left(\frac{b+ad-n/d}{c}\right)\epsilon\left(d,-\frac{b+ad-n/d}{c}\right)\right)q^n.
\end{equation}
Having in mind the above relations we get the proof of the theorem. $qed$\\
\\

Theorem 6 give us the Fourier coefficients of ($q=e(z)\textrm{, }Im(z)>0$):
$$
f_s(a,b;c;z)=
\sum^{\infty}_{n=-\infty}\frac{(-1)^nq^{an^2+bn}}{\sinh(2\pi i  n c z)}=\sum_{n\in\textbf{\scriptsize Z}-\{0\}}\frac{(-1)^ne\left((an^2+bn)z\right)}{\sinh(2\pi i n c z)}
$$
Lerch sums. In the same way the Lerch sum
\begin{equation}
f_c(a,b;c;z):=
\sum^{\infty}_{n=-\infty}\frac{(-1)^nq^{an^2+bn}}{\cosh(2\pi i  n c z)},
\end{equation}
can be written as:
$$
1+2\sum^{\infty}_{n=1}\frac{(-1)^nq^{an^2}\cosh(2\pi i b n  z)}{\cosh(2\pi i n c z)}=
$$
$$
1+2\sum^{\infty}_{l,n=1}(-1)^{l+n}q^{an^2+bn+cn+2lnc}+2\sum^{\infty}_{l,n=1}(-1)^{l+n}q^{an^2-bn+cn+2lnc}+
$$
$$
+2\sum^{\infty}_{n=1}(-1)^n q^{an^2+(c-b)n}+2\sum^{\infty}_{n=1}(-1)^n q^{an^2+(c+b)n}=
$$
$$
=1+2\sum^{\infty}_{l,n=-\infty}\epsilon_0(n,l)(-1)^{l+n}q^{an^2+c|n|+bn+2cnl}+
$$
$$
+2\sum^{\infty}_{n=1}(-1)^nq^{a n^2+(c-b)n}+2\sum^{\infty}_{n=1}(-1)^nq^{an^2+(c+b)n}=
$$
$$
1+2\sum^{\infty}_{l,n=-\infty}(-1)^{l+n}\epsilon_0(n,l)q^{an^2+c|n|+bn+2cnl}+
$$
$$
+2\sum^{\infty}_{n=1}(-1)^n q^{a n^2+(c-b)n}+2\sum^{\infty}_{n=1}(-1)^n q^{an^2+(c+b)n},
$$
where
\begin{equation}
\epsilon_0(n,l):=
\left\{\begin{array}{cc}
1\textrm{, if } nl>0\\
0\textrm{, else }
\end{array}\right\}.
\end{equation}
Proceeding with the same arguments, as we did with (53), we arrive to the next\\
\\
\textbf{Theorem 7.}\\
Assume that $a,b,c$ are integers with $a>0$ and $z$ is complex number with $Im(z)>0$. If $q=e(z)$, then the Lerch sum (54) have the following expansion in power series
$$
\sum^{\infty}_{n=-\infty}\frac{(-1)^nq^{an^2+bn}}{\cosh(2\pi i  n c z)}=1+\sum^{\infty}_{n=1}C_c(a,b;c;n)q^n+
2\sum^{\infty}_{n=1}(-1)^n q^{a n^2+(c-b)n}+
$$
\begin{equation}
+2\sum^{\infty}_{n=1}(-1)^n q^{an^2+(c+b)n},
\end{equation}   
where
\begin{equation}
C_c(a,b;c;n)=C_2(a,c+b;2c;n)+C_2(a,c-b;2c;n)
\end{equation}
and
$$
C_2(a,b;c;n):=
$$
\begin{equation}
2\sum_{d|n}(-1)^dX_{\textbf{\scriptsize{Z}}}\left(\frac{b+ad-n/d}{c}\right)\epsilon\left(d,-\frac{b+ad-n/d}{c}\right)(-1)^{-(b+ad-n/d)/c}.
\end{equation}
\\
\textbf{Remarks.}\\
1) Assume the theta function
\begin{equation}
\theta_{a,b}(\chi,q)=\sum^{\infty}_{n=1}\chi(n)q^{a n^2+b n}.
\end{equation}
We are going to evaluate its power series coefficients. For simplicity we assume $X(n,m)$ is such that $X(n,n)=\chi(n)$ and $\delta_{n,m}$ is 1, if $n=m$ and 0 else. Then 
$$
\theta_{a,b}(\chi,q)=\sum^{\infty}_{n,m=1}X(n,m)\delta_{n,m}q^{n(a m+b)}.
$$
Hence
\begin{equation}
\theta_{a,b}(\chi,q)=\sum^{\infty}_{n=1}q^n\sum_{d|n}X_{\textbf{\scriptsize{Z}}}\left(\frac{n/d-b}{a}\right)X\left(d,\frac{n/d-b}{a}\right)\delta_{d,(n/d-b)/a}.
\end{equation}
Consequently if $a,b$ are integers such $a,b>0$, or $a>|b|$, then
\begin{equation}
\theta_{a,b}(\chi,q)=\sum^{\infty}_{n=1}q^n\sum_{\scriptsize 
\begin{array}{cc}
	d|n\\
	ad^2+bd=n
\end{array}
\normalsize}\chi(d).
\end{equation} 
2) Therefore also
$$
C_c(a,b,c,n)=
$$
$$
2\sum_{\scriptsize
\begin{array}{cc}
	d|n\\
	n/d\equiv ad+c-b(2c)
\end{array}
\normalsize}(-1)^d\epsilon\left(d,-\frac{c-b+ad-n/d}{2c}\right)(-1)^{-(c-b+ad-n/d)/(2c)}+
$$
$$
+2\sum_{\scriptsize
\begin{array}{cc}
	d|n\\
	n/d\equiv ad+c+b(2c)
\end{array}
\normalsize}(-1)^d\epsilon\left(d,-\frac{c+b+ad-n/d}{2c}\right)(-1)^{-(c+b+ad-n/d)/(2c)}+
$$
\begin{equation}
+2\sum_{\scriptsize
\begin{array}{cc}
	d|n\\
	ad^2+(c-b)d=n
\end{array}\normalsize}(-1)^d+2\sum_{\scriptsize
\begin{array}{cc}
	d|n\\
	ad^2+(c+b)d=n
\end{array}
\normalsize}(-1)^d,
\end{equation}
where
\begin{equation}
\epsilon(n,l)=\textrm{sign}(n)\textrm{sign}(l)\frac{\textrm{sign}(n)+\textrm{sign}(l)}{2}
\end{equation}
and
\begin{equation}
\sum^{\infty}_{n=-\infty}\frac{(-1)^nq^{an^2+bn}}{\cosh(2\pi i  n c z)}=1+\sum^{\infty}_{n=1}C_c(a,b,c,n)q^n\textrm{, }q=e(z)\textrm{, }Im(z)>0.
\end{equation}

\section{More representations, evaluations and modular relations of Lerch sums}

In this section we will use Poisson summation formula to recover properties of Lerch sums. The Poisson summation formula read as (see [4]): 
\begin{equation}
\sum^{\infty}_{n=-\infty}f(n)=\sum^{\infty}_{n=-\infty}\widehat{f}(2\pi n),
\end{equation}
where
\begin{equation}
\widehat{f}(x)=f(.)\wedge(x)=\int^{+\infty}_{-\infty}f(t)e^{-itx}dt
\end{equation}
is the Fourier transform of $f(t)$.\\
Assume now $q:=e(z)=e^{2\pi i z}$. We define 
\begin{equation}
f_1(t)=\exp\left(2\pi i z (at^2+Bt)\right)\textrm{, }B=b-1
\end{equation}
and
\begin{equation}
f_2(t)=\frac{1}{\cosh(2\pi i t z^*)},
\end{equation}
where the asterisc ''$*$'' means complex conjugate. Then setting
\begin{equation}
S(a,b,z):=\sum^{\infty}_{n=-\infty}\frac{q^{an^2+bn}}{1+q^{2n}},
\end{equation}
we can write according to Poisson summation formula:
$$
S(a,b,z):=\frac{1}{2}\sum^{\infty}_{n=-\infty}\frac{e\left(z(an^2+(b-1)n)\right)}{\cosh(2 \pi i z n)}=\frac{1}{2}\sum^{\infty}_{n=-\infty}f_1(n)f_2^{*}(n)=
$$
\begin{equation}
=\frac{1}{2}\sum^{\infty}_{n=-\infty}\left(\int^{+\infty}_{-\infty}f_1(t)f_2^{*}(t)e^{-2\pi i n t}dt\right).
\end{equation}
If
\begin{equation}
F_1(n,t)=f_1(t)e^{-2\pi i n t}\textrm{ and }F_2(t)=f_2(t)
\end{equation}
then using Parseval identity we can write:
\begin{equation}
\xi(n,a,b,z):=\int^{+\infty}_{-\infty}f_1(t)e^{-2\pi i n t}f_2^{*}(t)dt=\frac{1}{2\pi}\int^{+\infty}_{-\infty}\widehat{F}_1(n,\gamma)\widehat{F}^{*}_2(\gamma)d\gamma.
\end{equation}
Since also
\begin{equation}
\widehat{F}_1(n,\gamma)=\frac{1}{\sqrt{-2 i a z}}\exp\left(\frac{-i\left(2\pi n+\gamma-2B\pi z\right)^2}{8a\pi z}\right)
\end{equation}
and
\begin{equation}
\widehat{F}^{*}_2(\gamma)=\frac{i}{2z}\cdot\textrm{sech}\left(\frac{i\gamma}{4z}\right),
\end{equation}
we have
$$
\xi(n,a,b,z)=\frac{1}{2\pi}\int^{+\infty}_{-\infty}\widehat{F}_1(n,\gamma)\widehat{F}^{*}_2(\gamma)d\gamma=
$$
\begin{equation}
=\frac{i}{4\pi z\sqrt{-2ia z}}\int^{+\infty}_{-\infty}\exp\left(-\frac{i(\gamma+2\pi(n-(b-1)z))^2}{8a\pi z}\right)\sec\left(\frac{\gamma}{4z}\right)d\gamma.
\end{equation}
Hence (70) can be written as
$$
S(a,b,z)=\frac{1}{4\pi}\sum^{\infty}_{n=-\infty}\int^{+\infty}_{-\infty}\widehat{F}_1(n,\gamma)\widehat{F}^*_2(\gamma)d\gamma=
$$
\begin{equation}
\frac{1}{2}\sum^{\infty}_{n=-\infty}\frac{i}{4\pi z\sqrt{-2ia z}}\int^{+\infty}_{-\infty}\exp\left(-\frac{i(\gamma+2\pi(n-(b-1)z))^2}{8a\pi z}\right)\sec\left(\frac{\gamma}{4z}\right)d\gamma
\end{equation}
We can make also the following change of variable $\gamma=h L$ in (75) to get
$$
\frac{iL}{4\pi z\sqrt{-2iaz}}\int^{+\infty/L}_{-\infty/L}\exp\left(-\frac{i(hL+2\pi(n+(b-1)z))^2}{8a\pi z}\right)\frac{dh}{\cos\left(\frac{hL}{4z}\right)}=
$$
$$
=\frac{i L}{4\pi z\sqrt{-2iaz}}\exp\left(-\frac{i\pi(n-(b-1)z)^2}{2az}\right)\times
$$
\begin{equation}
\times\int^{+\infty/L}_{-\infty/L}\exp\left[-\frac{iL^2}{8a\pi z}h^2-\frac{iL}{2az} \left(n-z(b-1)\right)h\right]\frac{dh}{\cosh(i h L/(4z))}.
\end{equation}
But also
$$
\xi\left(n,a,b,z\right)=
$$
\begin{equation}
=\int^{+\infty}_{-\infty}\exp\left(2\pi i z a t^2-2\pi i  \left(n-z(b-1)\right)t\right)\frac{dt}{\cosh\left(2\pi i z t\right)}.
\end{equation}

Hence we have the next definition-theorem:\\ 
\\
\textbf{Theorem 8.}\\
If $Re(2\pi i z a)<0$, then 
\begin{equation}
\xi\left(n,a,b,z\right)=\int^{+\infty}_{-\infty}\exp\left(2\pi i z a t^2-2\pi i  \left(n-(b-1)z\right)t\right)\frac{dt}{\cosh\left(2\pi i z t\right)}
\end{equation}
and if $Re\left(-\frac{iL^2}{8a \pi z}\right)<0$, then
$$
\xi(n,a,b,z)
=\frac{i L}{4\pi z\sqrt{-2iaz}}\exp\left(-\frac{i\pi(n-(b-1)z)^2}{2az}\right)\times
$$
\begin{equation}
\times\int^{+\infty/L}_{-\infty/L}\exp\left[-\frac{iL^2}{8a\pi z}h^2-\frac{iL}{2az} \left(n-(b-1)z\right)h\right]\frac{dh}{\cosh(i h L/(4z))}.
\end{equation}
\\

Now we rearange the order of summation and integration in (76) to express $S(a,b,z)$ with a simple integral. For $Im(w)>0$ it holds
\begin{equation}
\sum^{\infty}_{n=-\infty}\exp\left(-\frac{i(t+2n\pi)^2}{8\pi w}\right)=\sqrt{-2iw}\cdot\vartheta_3\left(t/2,e^{2\pi i w}\right), 
\end{equation}
where
\begin{equation}
\vartheta_3(z,q):=\sum^{\infty}_{n=-\infty}q^{n^2}e^{2niz}\textrm{, }|q|<1.
\end{equation}
Hence using (81) in (76), we get
$$
S(a,b;z)=\frac{1}{2}\sum^{\infty}_{n=-\infty}\xi(n,a,b,z)=
$$
$$
=\frac{i}{4}\int^{+\infty/z}_{-\infty/ z}\vartheta_3\left((h+1-b)\pi z,e(az)\right)\sec\left(\frac{h\pi}{2}\right)dh=
$$
$$
=\frac{i}{4}\int^{+\infty/z}_{-\infty/ z}\vartheta_3\left(h\pi z,e(az)\right)\sec\left(\frac{(h+b-1)\pi}{2}\right)dh=
$$
\begin{equation}
=\frac{i}{4\pi z}\int^{+\infty}_{-\infty}\vartheta_3\left(h,e(az)\right)\sec\left(\frac{h}{2z}+\frac{(b-1)\pi}{2}\right)dh.
\end{equation}
The changing of sum and integration needs explanation. One can see that, for $a>0$ and $Im(z)>0$ the series 
$$
\sum^{\infty}_{n=-\infty}\exp\left(-\frac{i(\gamma+2\pi n-2\pi(b-1) z)^2}{8a\pi z}\right)
$$
is uniformly convergent (Weierstrass test). This justifies the rearanging of summation and integration.\\ 
Also it is known (see [1]) that for the transformation of variables  
\begin{equation}
a'=\frac{1}{a}\textrm{, }z'=\frac{-1}{4z}\textrm{, }w'=2wa'z',
\end{equation}
holds
\begin{equation}
\vartheta_3\left(w',e\left(a'z'\right)\right)=\sqrt{-2iaz}\exp\left(\frac{i w^2}{2\pi a z}\right)\vartheta_3\left(w,e(az)\right).
\end{equation}
Hence\\
\\ 
\textbf{Theorem 9.}\\
Setting
\begin{equation}
p_j(a,b,z):=\frac{i}{4\pi z}\int^{+\infty a}_{-\infty a}\vartheta_3\left(h,e(az)\right)\exp\left(j\frac{ih^2}{2\pi a z}\right)\sec\left(\frac{h}{2z}+\frac{b-1}{2}\pi\right)dh,
\end{equation}
where $j=0,1$. Then
\begin{equation}
p_0\left(\frac{1}{2a},b,\frac{-1}{2z}\right)=-\sqrt{-2iaz}\cdot p_1\left(z,b,a\right)\textrm{, where }a,Im(z)>0
\end{equation}
and also
\begin{equation}
p_0(a,b,z)=\frac{1}{2}\sum^{\infty}_{n=-\infty}\frac{e\left((an^2+(b-1)n)z\right)}{\cosh(2\pi i n z)}=\sum^{\infty}_{n=-\infty}\frac{q^{an^2+bn}}{1+q^{2n}}.
\end{equation}
\\
\textbf{Proof.}\\
Relation (88) is consequence of the evaluations (83). The proof of (87) is a straight forward application of the transformations (84) and (85) of $\vartheta_3(h,e(az))$ functions. $qed$\\ 
\\
\textbf{Note.}\\
I could not find any physical meaning or where it might leads relation (87). It seems to me very weird.\\
\\
\textbf{Theorem 10.}\\
In the same way as above if $a>0$, $Im(z)>0$, $Im(w)\neq0$, $q=e(z)$, we have:
\begin{equation}
\sum^{\infty}_{n=-\infty}\frac{q^{an^2+bn}}{\cosh(2\pi i n w)}=\frac{i}{2\pi w}\int^{+\infty}_{-\infty }\vartheta_3\left(h,e(az)\right)\sec\left(\frac{h+b\pi z}{2w}\right)dh
\end{equation}
and
\begin{equation}
\sum^{\infty}_{n=-\infty}\frac{(-1)^nq^{an^2+bn}}{\cosh(2\pi i n w)}=\frac{i}{2\pi w}\int^{+\infty}_{-\infty }\vartheta_4\left(h,e(az)\right)\sec\left(\frac{h+b\pi z}{2w}\right)dh.
\end{equation}
\\

We proceed now with a definition. For $j=0,1$, we define more general
$$
P_j(a,b;z,w;x):=
$$
\begin{equation}
=\frac{i}{2\pi w}\int^{+\infty a}_{-\infty a}\vartheta_3\left(t,e(az)\right)\exp\left(j\frac{it^2}{2\pi a z}\right)\sec\left(\frac{t}{2w}+\frac{b\pi z}{2w}\right)e^{it x}dt.
\end{equation}
From Theorem 10 is
\begin{equation}
P_0(a,b;z,w;0)=\sum^{\infty}_{n=-\infty}\frac{q^{an^2+b n}}{\cosh(2\pi i n w)}\textrm{, }a>0.
\end{equation}
This function is a Lerch series function. Hence we can say that $P_j(a,b;z,w;x)$ defined above is a generalization of Lerch series. The Fourier transformation of these generalized functions possess modular properties. Moreover we can prove the next:\\
\\
\textbf{Theorem 11.}\\
If $a>0$, $Im(z)>0$ and $Im(w)\neq0$, then
\begin{equation}
P_{1-j}(a',b';z',w';(.))\wedge(\gamma')
=-i (-2ia z)^{3/2}P_{j}(a,b;z,w;(.))\wedge(\gamma),
\end{equation}
where $j=0,1$ and
\begin{equation}
a'=\frac{1}{a}\textrm{, }b'=2ba'z\textrm{, }z'=\frac{-1}{4z}\textrm{, }\gamma'=2\gamma a' z'\textrm{, }w'=2wa'z'.
\end{equation}
\\
\textbf{Proof.}\\
When $a>0$ and $j=0,1$, from relation (91) we have
\begin{equation}
P_j\left(a,b;z,w;(.)\right)\wedge(\gamma)=\frac{i}{ w}\vartheta_3\left(\gamma,e(az)\right)\exp\left(j\frac{i\gamma^2}{2\pi a z}\right)\sec\left(\frac{\gamma}{2w}+\frac{b\pi z}{2w}\right).
\end{equation}
Using the transformation property of $\vartheta_3(t,e(az))$ (relation (85)) and relations (94), we get the result. $qed$\\
\\

Assume now that $a,Im(z)>0$, $Im(w)\neq0$, then 
$$
P_{1-j}(a',b';z',w';(.))\wedge(0)=-2a\sqrt{-2ia}\cdot z^{3/2}P_j(a,b;z,w;(.))\wedge(0).
$$
Setting $\gamma=0$ in relation (95) we have
$$
P_j(a,b;z,w;(.))\wedge(0)=\frac{i}{ w}\vartheta_3(0,e(az))\sec\left(\frac{b\pi z}{2w}\right). 
$$
From (94),(95) we have 
$$
P_{1}(a,b;z,w;x)=\frac{1}{2\pi}\int^{+\infty a}_{-\infty a}P_{0}(a,b;z,w;(.))\wedge(\gamma)\exp\left(\frac{i\gamma^2}{2\pi a z}\right)e^{i\gamma x}d\gamma.
$$
and
$$
P_{0}(a,b;z,w;x)=\frac{1}{2\pi}\int^{+\infty a}_{-\infty a}P_{1}(a,b;z,w;(.))\wedge(\gamma)\exp\left(-\frac{i\gamma^2}{2\pi a z}\right)e^{i\gamma x}d\gamma.
$$
Hence
\begin{equation}
P_{1-j}(a,b;z,w;x)=\frac{1}{2\pi}\int^{+\infty a}_{-\infty a}P_{j}(a,b;z,w;(.))\wedge(\gamma)\exp\left((-1)^j\frac{i\gamma^2}{2\pi a z}\right)e^{i\gamma x}d\gamma.
\end{equation}
Hence if $a>0$ using the isometry property of the Fourier transform in a convolution product, we get for $j=0,1$:
\begin{equation}
P_{1-j}(a,b;z,w;x)
=\sqrt{\frac{iaz }{2(-1)^j}}\int^{+\infty}_{-\infty}P_j\left(a,b;z,w;x-t\right)\exp\left(\frac{-i t^2 \pi a z}{2(-1)^j}\right)dt.
\end{equation}
Also
$$
P_{1-j}(a',b';z',w';x)=
$$
$$
=\frac{i\sqrt{-2iaz}}{2\pi w}\int^{+\infty z}_{-\infty z}\vartheta_3\left(t,e(az)\right)\exp\left(j\frac{it^2}{2\pi a z}\right)\sec\left(\frac{t}{2w}+\frac{b\pi z}{2w}\right)e^{-itx/(2az)}dt.
$$
Hence
\begin{equation}
P_{1-j}(a',b';z',w';x)=\sqrt{-2ia z}\cdot P_{j}\left(a,b;z,w;\frac{-x}{2az}\right)
\end{equation}
and equivalently
\begin{equation}
P_{j}(a',b';z',w';-2xaz)=\sqrt{-2ia z}\cdot P_{1-j}\left(a,b;z,w;x\right).
\end{equation}
Consequently from (97):\\
\\
\textbf{Theorem 12.}\\
If $a>0$, $Im(z)>0$, $Im(w)\neq 0$ and 
\begin{equation}
a'=\frac{1}{a}\textrm{, }b'=2ba'z'\textrm{, }z'=\frac{-1}{4z}\textrm{, }w'=2wa'z'\textrm{, }x'=-2xaz,
\end{equation}
then for $j=0,1$ we have
\begin{equation}
P_{j}(a',b';z',w';x')=\frac{az }{i^j}\int^{+\infty}_{-\infty}P_j\left(a,b;z,w;x-t\right)\exp\left(\frac{-i t^2 \pi a z}{2(-1)^j}\right)dt.
\end{equation} 
\\
\textbf{Examples.}\\
1) Suppose that $a=3/2$, $b=0$, $w=z/2$, then
$$
f(q)=\frac{2}{\eta(q)}\sum^{\infty}_{n=-\infty}\frac{(-1)^n q^{3n^2/2+n/2}}{1+q^n}=
$$
\begin{equation}
=\frac{i}{\pi z\cdot \eta\left(e(z)\right)}\int^{+\infty}_{-\infty}\vartheta_4\left(h,e\left(\frac{3z}{2}\right)\right)\sec\left(\frac{h}{z}\right)dh.
\end{equation}
\\
2) Let $q=e(z)$, $Im(z)>0$, then
$$
\eta(q)\phi(q)=\eta(q)\sum^{\infty}_{n=0}\frac{q^{n^2}}{\left(-q^2;q^2\right)_n}=\sum^{\infty}_{n=-\infty}\frac{(-1)^n(1+q^n)q^{n(3n+1)/2}}{1+q^{2n}}=
$$
$$
=\sum^{\infty}_{n=-\infty}\frac{(-1)^nq^{n(3n+1)/2}}{1+q^{2n}}+\sum^{\infty}_{n=-\infty}\frac{(-1)^nq^{3n(n+1)/2}}{1+q^{2n}}.
$$
Hence
$$
2\eta\left(q^2\right)\phi\left(q^2\right)=\sum^{\infty}_{n=-\infty}\frac{(-1)^nq^{3n^2-n}}{\cosh(4\pi i n z)}+\sum^{\infty}_{n=-\infty}\frac{(-1)^n q^{3n^2+n}}{\cosh(4\pi i n z)}
$$
Observe now that the two sums are the same and we can write
\begin{equation}
\eta\left(q^2\right)\phi\left(q^2\right)=\sum^{\infty}_{n=-\infty}\frac{(-1)^nq^{3n^2+n}}{\cosh(4\pi i nz)}.
\end{equation}
For evaluating the Fourier coefficients of the above Lerch sum we have
$$
\eta\left(q^2\right)\phi\left(q^2\right)=
$$
$$
=1+\sum^{\infty}_{n=1}C_c(3,1;2;n)q^n+2\sum^{\infty}_{n=1}(-1)^nq^{3n(n+1)}+2\sum^{\infty}_{n=1}(-1)^nq^{3n^2+n},
$$
where the arithmetic function $C_c$ is given in Theorem 7. Using also (60),(61), we get
$$
\eta\left(q^2\right)\phi\left(q^2\right)=1+\sum^{\infty}_{n=1}C_c(3,1;2;n)q^n+
$$
$$
+2\sum^{\infty}_{n=1}\left(\sum_{d|n}X_{\textbf{\scriptsize{Z}}}\left(\frac{n/d-3}{3}\right)X\left(d,\frac{n/d-3}{3}\right)\delta_{d,(n/d-3)/3}\right)q^n+
$$
$$
+2\sum^{\infty}_{n=1}\left(\sum_{d|n}X_{\textbf{\scriptsize{Z}}}\left(\frac{n/d-1}{3}\right)X\left(d,\frac{n/d-1}{3}\right)\delta_{d,(n/d-1)/3}\right)q^n,
$$
where $X(n,m)$ is any arithmetic function such that $X(n,n)=(-1)^n$. An example is $X(n,m)=(-1)^{(n+m)/2}$. Hence
$$
\eta\left(q^2\right)\phi\left(q^2\right)=
$$
$$
=1+2\sum^{\infty}_{n=1}[\sum_{\scriptsize
\begin{array}{cc}
	d|n\\
	n/d\equiv 3d+1(4)
\end{array}
\normalsize}(-1)^d\epsilon\left(d,-\frac{1+3d-n/d}{4}\right)(-1)^{-(1+3d-n/d)/4}+
$$
$$
+\sum_{\scriptsize
\begin{array}{cc}
	d|n\\
	n/d\equiv 3d+3(4)
\end{array}
\normalsize}(-1)^d\epsilon\left(d,-\frac{3+3d-n/d}{4}\right)(-1)^{-(3+3d-n/d)/4}+
$$
\begin{equation}
+\sum_{\scriptsize
\begin{array}{cc}
	d|n\\
	3d^2+d=n
\end{array}\normalsize}(-1)^d+\sum_{\scriptsize
\begin{array}{cc}
	d|n\\
	3d^2+3d=n
\end{array}
\normalsize}(-1)^d] q^n.
\end{equation}
 
Continuing we can write (the integral representation):
$$
2\eta(q)\phi(q)=\int^{+\infty}_{-\infty}\vartheta_4\left(h,e\left(\frac{3z}{2}\right)\right)\sec\left(\frac{h}{2z}-\frac{\pi}{4}\right)dh+
$$
$$
+\int^{+\infty}_{-\infty}\vartheta_4\left(h,e\left(\frac{3z}{2}\right)\right)\sec\left(\frac{h}{2z}+\frac{\pi}{4}\right)dh.
$$
Hence
\begin{equation}
\phi(q)=\frac{i\sqrt{2}}{2\pi z\cdot  \eta(q)}\int^{+\infty}_{-\infty}\vartheta_4\left(h,e\left(\frac{3z}{2}\right)\right)\cos\left(\frac{h}{2z}\right)\sec\left(\frac{h}{z}\right)dh.
\end{equation}

\section{The general case of representation of Lerch sums}

From [3] we have the next evaluation formula
\begin{equation}
\sum^{\infty}_{n=-\infty}(-1)^nq^{pn^2/2+(p-2a)n/2}=q^{-p/12+a/2-a^2/(2p)}\eta(q^p)W^{\{4\}}_{\{a,p\}}(m(q)),
\end{equation}
where
$$
W^{\{4\}}_{\{a,p\}}\left(m(q)\right)=A(a,p;q)=
$$
$$
=q^{p/12-a/2+a^2/(2p)}\prod^{\infty}_{n=0}(1-q^{np+a})(1-q^{np+p-a})=
$$
\begin{equation}
=q^C[a,p;q]^{+}_{\infty}=q^C\left(q^a;q^p\right)_{\infty}\left(q^{p-a};q^p\right)_{\infty}.
\end{equation}
But 
\begin{equation}
\log\left(\left(q^a;q^p\right)_{\infty}\left(q^{p-a};q^p\right)_{\infty}\right)=-\sum^{\infty}_{n=1}q^n\sum_{\scriptsize
\begin{array}{cc}
AB=n\\
B\equiv \pm a(p)	
\end{array}}\frac{1}{A}.
\end{equation}
Hence
\begin{equation}
W^{\{4\}}_{\{a,p\}}\left(m(q)\right)=q^{p/12-a/2+a^2/(2p)}\exp\left(-\sum^{\infty}_{n=1}q^n\sum_{\scriptsize
\begin{array}{cc}
AB=n\\
B\equiv \pm a(p)	
\end{array}\normalsize}\frac{1}{A}\right).
\end{equation}
Also in the same way if $|q|<1$, then
\begin{equation}
\sum^{\infty}_{n=-\infty}q^{pn^2/2+(p-2a)n/2}=q^{-p/12+a/2-a^2/(2p)}\eta(q^p)W^{\{3\}}_{\{a,p\}}(m(q)),
\end{equation}
where
\begin{equation}
W^{\{3\}}_{\{a,p\}}\left(m(q)\right)=q^{p/12-a/2+a^2/(2p)}\exp\left(-\sum^{\infty}_{n=1}q^n\sum_{\scriptsize
\begin{array}{cc}
AB=n\\
B\equiv \pm a(p)	
\end{array}\normalsize}\frac{(-1)^A}{A}\right).
\end{equation}
We want to find the analog of (106) for the Jacobi theta function $\vartheta_4$. This can be done considering the next transformation of variables: $a\rightarrow a-t$ and $p\rightarrow 2a$ in (106). Then
$$
\vartheta_4\left(\pi zt,e(az)\right)=q^{a/12-t^2/(4a)}\eta\left(q^{2a}\right)Q^{\{4\}}_{\{a,t\}}\left(m(q)\right).
$$
Hence if $q=e(z)$, $Im(z)>0$, then
\begin{equation}
\vartheta_4\left(\pi z t ,e(az)\right)=q^{-t^2/(4a)}\eta_D\left(2az\right)Q^{\{4\}}_{\{a,t\}}\left(m(q)\right),
\end{equation}
where
\begin{equation}
Q^{\{4\}}_{\{a,t\}}(m(q))=q^{-a/12+t^2/(4a)}\prod^{\infty}_{n=0}\left(1-q^{2na+a-t}\right)\left(1-q^{2na+a+t}\right).
\end{equation}
In case of $a,t\in\textbf{Z}$, $a>t>0$, then
\begin{equation}
Q^{\{4\}}_{\{a,t\}}(m(q))=q^{-a/12+t^2/(4a)}\exp\left(-\sum^{\infty}_{n=1}q^n\sum_{\scriptsize
\begin{array}{cc}
AB=n\\
B\equiv \pm (a-t)(mod 2a)	
\end{array}\normalsize}\frac{1}{A}\right)
\end{equation}  
and
$$
m(q):=\left(\frac{\vartheta_2(0,q)}{\vartheta_3(0,q)}\right)^2\textrm{, }|q|<1.
$$
Also if $q=e(z)$, $Im(z)>0$, then
\begin{equation}
\vartheta_3\left(\pi zt,e(az)\right)=q^{-t^2/(4a)}\eta_D\left(2az\right)Q^{\{3\}}_{\{a,t\}}\left(m(q)\right),
\end{equation}
where
\begin{equation}
Q^{\{3\}}_{\{a,t\}}(m(q))=q^{-a/12+t^2/(4a)}\prod^{\infty}_{n=0}\left(1+q^{2na+a-t}\right)\left(1+q^{2na+a+t}\right).
\end{equation}
In case $a,t\in\textbf{Z}$, $a>t>0$, then
\begin{equation}
Q^{\{3\}}_{\{a,t\}}(m(q))=q^{-a/12+t^2/(4a)}\exp\left(-\sum^{\infty}_{n=1}q^n\sum_{\scriptsize
\begin{array}{cc}
AB=n\\
B\equiv \pm (a-t)(mod 2a)	
\end{array}\normalsize}\frac{(-1)^A}{A}\right).
\end{equation} 
But it holds the following modular identity
\begin{equation}
\vartheta_3\left(\pi t' z',e(a'z')\right)=\sqrt{-2iaz}\exp\left(\frac{i\pi t^2 z}{2a}\right)\vartheta_3\left(\pi t z,e(az)\right),
\end{equation}
where
\begin{equation}
a'=1/a\textrm{, }z'=-1/(4z)\textrm{, }t'=2tz/a.
\end{equation}
Hence in general for the function $F_3(a,t;z):=Q^{\{3\}}_{\{a,t\}}(m(q))$, $q=e(z)$ it holds
$$
\frac{F_3(a',t';z')}{F_3(a,t;z)}=\sqrt{-2iaz}\exp\left(\frac{-i\pi t^2 z}{2a}\right)\frac{\eta_D(2az)}{\eta_D\left(\frac{-1}{2az}\right)},
$$
where $\eta_D(z)$, $Im(z)>0$ is the Dedekind's eta function. Using the next functional equation: 
\begin{equation}
\eta_D\left(-1/z\right)=\sqrt{-iz}\cdot\eta_D(z),
\end{equation}
we finally arrive to\\
\\
\textbf{Theorem 13.}\\
1) $(Conjecture)$ If $a>0$ and $Im(z)>0$ 
\begin{equation}
\vartheta_3\left(\pi t z,e(az)\right)=q^{a/12-t^2/(4a)}\eta\left(q^{2a}\right)F_3(a,t;z),
\end{equation}
the function $F_3(a,t;z)$ takes algebraic values when $a,t\in\textbf{Q}^{*}_{+}$ and $z=r_1+i\sqrt{r_2}$, with $r_1$ rational and $r_2$ is positive rational.\\
2) If $a,t$ positive integers with $a>t$, then
$$
F_3(a,t;z)=Q^{\{3\}}_{\{a,t\}}(m(q))=
$$
\begin{equation}
=q^{-a/12+t^2/(4a)}\exp\left(-\sum^{\infty}_{n=1}q^n\sum_{\scriptsize
\begin{array}{cc}
AB=n\\
B\equiv \pm (a-t)(mod 2a)	
\end{array}\normalsize}\frac{(-1)^A}{A}\right).
\end{equation}
Relation (122) is valid only for $a,t$ positive integers. However we can use other ways to work, if the problem requires it. Such formulas are (116) or the definition (121), which are valid in reals.\\ 
3) For the transformation of variables (119) it holds
\begin{equation}
F_3(a',t';z')=\exp\left(\frac{-i\pi t^2 z}{2a}\right)F_3(a,t;z).
\end{equation}
\\
\textbf{Theorem 14.} $(Conjecture)$\\
When $a>t$ and $a,t$ positive rationals, the function $Q^{\{y\}}_{\{a,t\}}(x)$, $y=3,4$ takes algebraic numbers to algebraic numbers.\\
\\
\textbf{Theorem 15.}\\
If $q=e(z)$ and $a>0$, $Im(z)>0$, $Im(w)\neq0$, then
$$
S=\frac{1}{\eta_D\left(2az\right)}\sum^{\infty}_{n=-\infty}\frac{q^{an^2+bn}}{\cosh(2\pi i n w)}=
$$
\begin{equation}
=\frac{i}{2w}\int^{+\infty}_{-\infty}F_{3}\left(\frac{1}{a},\frac{2h}{a};\frac{-1}{4z}\right)\sec\left(\frac{(h+bz)\pi}{2w}\right)dh,
\end{equation}
where $\eta_D(z)$ is the classical Dedekind eta function.\\
\\
\textbf{Proof.}\\
From Theorems 10 and 13 we can write
$$
\frac{1}{\eta_D\left(2az\right)}\sum^{\infty}_{n=-\infty}\frac{q^{an^2+bn}}{\cosh(2\pi i n w)}=
$$
$$
=\frac{i\pi z}{2\pi w}\int^{+\infty/z}_{-\infty/z}\frac{\vartheta_3(\pi h z,e(az))}{\eta_D\left(2az\right)}\sec\left(\frac{(h+b)\pi z}{2w}\right)dh=
$$
$$
=\frac{iz}{2 w}\int^{+\infty/z}_{-\infty/z}F_3(a,h,z)e^{-\pi i h^2 z/(2a)}\sec\left(\frac{(h+b)\pi z}{2w}\right)dh
=
$$
$$
=\frac{i}{2 w}\int^{+\infty}_{-\infty}F_3\left(a,\frac{h}{z};z\right)e^{-i\pi h^2/(2az)}\sec\left(\frac{(h\pi+b\pi z}{2w}\right)dh=
$$
$$
=\frac{i}{2 w}\int^{+\infty}_{-\infty}F_3\left(\frac{1}{a},\frac{2h}{a},\frac{-1}{4z}\right)\sec\left(\frac{(h\pi+b\pi z}{2w}\right)dh,
$$
where we have used modular relation (123) to arive to the desired result. $qed$\\
\\
\textbf{Theorem 16.}\\
For $j=0,1$ we define
$$
S_j(a,b;z,w;\sigma)=\frac{i}{2w}\int^{+\infty \sigma}_{-\infty \sigma}F_3\left(\frac{1}{a},\frac{2h}{a};\frac{-1}{4z}\right)e^{i \pi h^2 j/(2az)}\sec\left(\frac{(h+b z)\pi}{2w}\right)dh,
$$
then for the transformation of variables
\begin{equation}
a'=1/a\textrm{, }b'=2ba'z\textrm{, }z'=-1/(4z)\textrm{, }w'=2wa'z',
\end{equation}
holds
\begin{equation}
S_1\left(a',b';z',w';z'\right)=S_0\left(a,b;z,w;a\right).
\end{equation}
Also when $a>0$, then
\begin{equation}
S_0(a,b;z,w;a)=\frac{1}{\eta_D\left(2az\right)}\sum^{\infty}_{n=-\infty}\frac{q^{an^2+bn}}{\cosh(2\pi i n w)}.
\end{equation}
\\
\textbf{Proof.}\\
$$
S_1\left(a',b';z',w';z'\right)=
$$
$$
=\frac{i}{2w'}\int^{+\infty z'}_{-\infty z'}F_3\left(\frac{1}{a'},\frac{2h}{a'};\frac{-1}{4z'}\right)e^{ih^2\pi/(2a'z')}\sec\left(\frac{h\pi}{2w'}+\frac{b'\pi z'}{2w'}\right)dh=
$$
$$
=\frac{i}{2w'}\int^{-\infty/z}_{+\infty/z}F_3\left(a,2ha;z\right)e^{-2ih^2az\pi}\sec\left(\frac{h\pi}{2w'}-\frac{b'\pi}{8w'z}\right)dh=
$$
$$
=\frac{i}{4w'a}\int^{-\infty a/z}_{+\infty a/z}F_3\left(a,h;z\right)e^{-ih^2 z \pi/(2a)}\sec\left(\frac{h\pi}{4aw'}-\frac{b'\pi}{8zw'}\right)dh=
$$
$$
=\frac{i}{4w'a}\int^{-\infty a/z}_{+\infty a/z}F_3\left(\frac{1}{a},\frac{2hz}{a};\frac{-1}{4z}\right)\sec\left(\frac{h\pi }{4aw'}-\frac{b'\pi}{8zw'}\right)dh=
$$
$$
=-\frac{i}{4w'az}\int^{+\infty a}_{-\infty a}F_3\left(\frac{1}{a},\frac{2h}{a};\frac{-1}{4z}\right)\sec\left(\frac{h\pi}{4aw'z}-\frac{b'\pi}{8zw'}\right)dh=
$$
$$
=\frac{i}{2w}\int^{+\infty a}_{-\infty a}F_3\left(\frac{1}{a},\frac{2h}{a};\frac{-1}{4z}\right)\sec\left(\frac{h\pi}{2w}+\frac{b\pi z}{2w}\right)dh. 
$$
$qed$\\
\\

In the same way as above we can consider the Fourier transform pairs
\begin{equation}
\exp\left(2\pi i z(a(.)^2+b(.))-2\pi i n(.)\right)\leftrightarrow\frac{1}{\sqrt{-2ia z}}\exp\left(\frac{-i(2n\pi+(.)-2b\pi z)^2}{8a\pi z}\right)
\end{equation}
and
\begin{equation}
\frac{1}{\cosh(2\pi i w (A(.)+B))}\leftrightarrow i\frac{\exp\left(\frac{iB(.)}{A}\right)}{2Aw}\sec\left(\frac{(.)}{4Aw}\right),
\end{equation}
with $a,A,B,Im(z)>0$, $Im(w)\neq0$. Then we can write ($q=e(z)$):
$$
S_1=\sum^{\infty}_{n=-\infty}\frac{q^{an^2+bn}}{\cosh(2\pi i w^{*}(An+B))}=\sum^{\infty}_{n=-\infty}f_1(n)f^{*}_2(n),
$$
where
$$
f_1(t)=q^{at^2+bt}\textrm{ and }f_2(t)=\frac{1}{\cosh(2\pi i w(At+B))}.
$$
Hence
$$
S_1=\sum^{\infty}_{n=-\infty}\int^{+\infty}_{-\infty}f_1(t)e^{-2\pi i n t}f_2^{*}(t)dt
=\sum^{\infty}_{n=-\infty}\int^{+\infty}_{-\infty}F_1(n,t)F_2^{*}(t)dt,
$$
where $F_1(n,t)=e(z(at^2+bt)-nt)$ and $F_2(t)=f_2(t)$. Hence from the Parseval's identity and relations (128),(129), we get
$$
S_1=\frac{1}{2\pi}\sum^{\infty}_{n=-\infty}\int^{+\infty}_{-\infty}\widehat{F}_1(n,\gamma)\widehat{F}_2^{*}(\gamma)d\gamma=
$$
$$
=\frac{i}{4\pi Aw^{*} \sqrt{-2ia z}}\sum^{\infty}_{n=-\infty}\int^{+\infty}_{-\infty}e^{-i(2n\pi+\gamma-2b\pi z)^2/(8a\pi z)}
e^{-iB\gamma/A}\sec\left(\frac{\gamma}{4Aw^{*}}\right)d\gamma.
$$
Rearanging the order of summation and integration and using the formula
$$
\frac{1}{\sqrt{-2ia z}}\sum^{\infty}_{n=-\infty}\exp\left(\frac{-i(2n\pi+\gamma-2b\pi z)^2}{8a\pi z}\right)=\vartheta_3\left(\frac{\gamma}{2}-b\pi z,e(az)\right),
$$
we get
$$
S_1=\frac{-i}{4A\pi w^{*}}\int^{+\infty}_{-\infty}\vartheta_3\left(\frac{\gamma}{2}-b\pi z,e(az)\right)e^{-iB\gamma/A}\sec\left(\frac{\gamma}{4Aw^{*}}\right)d\gamma=
$$
$$
=\frac{1}{2A\pi i w^{*}}\int^{+\infty}_{-\infty}\vartheta_3\left(\gamma-b\pi z,e(az)\right)e^{-2iB\gamma/A}\sec\left(\frac{\gamma}{2Aw^{*}}\right)d\gamma=
$$
$$
=\frac{b\pi z}{2A\pi i w^*}\int^{+\infty/z}_{-\infty/z}\vartheta_3\left((\gamma-1)b\pi z,e(az)\right)e^{-2iB\gamma b\pi z/A}\sec\left(\frac{\gamma b\pi z}{2Aw^*}\right)d\gamma=
$$
$$
=\frac{b\pi z}{2A\pi i w^*}\int^{+\infty/z}_{-\infty/z}\vartheta_3\left(\gamma b\pi z,e(az)\right)e^{-2iB(\gamma+1)b\pi z/A}\sec\left(\frac{(\gamma+1) b\pi z}{2Aw^*}\right)d\gamma=
$$
$$
=\frac{b\pi z}{2A\pi i w^*}e^{-2iBb\pi z/A}\int^{+\infty/z}_{-\infty/z}\vartheta_3\left(\gamma b\pi z,e(az)\right)e^{-2iB\gamma b\pi z/A}\sec\left(\frac{(\gamma+1) b\pi z}{2Aw^*}\right)d\gamma=
$$
$$
=\frac{e^{-2iBb\pi z/A}}{2A\pi i w^*}\int^{+\infty}_{-\infty}\vartheta_3\left(\gamma,e(az)\right)e^{-2iB\gamma/A}\sec\left(\frac{\gamma+b\pi z}{2Aw^*}\right)d\gamma.
$$
Hence we have proven the next\\
\\
\textbf{Theorem 17.}\\
If $q=e(z)$, $Im(z)>0$, $a>0$ and $Im(w)\neq 0$, then
$$
\sum^{\infty}_{n=-\infty}\frac{q^{an^2+bn}}{\cosh(2\pi i w(An+B))}=
$$
\begin{equation}
=\frac{e^{-2iBb\pi z/A}}{2A\pi i w}\int^{+\infty}_{-\infty}\vartheta_3\left(\gamma,e(az)\right)e^{-2iB\gamma/A}\sec\left(\frac{\gamma+b\pi z}{2Aw}\right)d\gamma
\end{equation}
and
$$
\sum^{\infty}_{n=-\infty}\frac{(-1)^nq^{an^2+bn}}{\cosh(2\pi i w(An+B))}=
$$
\begin{equation}
=\frac{e^{-2iBb\pi z/A}}{2A\pi i w}\int^{+\infty}_{-\infty}\vartheta_4\left(\gamma,e(az)\right)e^{-2iB\gamma/A}\sec\left(\frac{\gamma+b\pi z}{2Aw}\right)d\gamma.
\end{equation}
Also as in Theorem 15 we get
$$
\frac{1}{\eta_D(2az)}\sum^{\infty}_{n=-\infty}\frac{q^{an^2+bn}}{\cosh(2\pi i w(An+B))}=
$$
\begin{equation}
=\frac{e^{-2iBb\pi z/A}}{2A i w}\int^{+\infty}_{-\infty}F_3\left(\frac{1}{a},\frac{2\gamma}{a};\frac{-1}{4z}\right)e^{-2i\pi B\gamma/A}\sec\left(\frac{(\gamma+b z)\pi}{2Aw}\right)d\gamma.
\end{equation}
\\
\textbf{Example.}\\
Assume the Ramanujan mock theta function 
\begin{equation}
\psi(q)=\sum^{\infty}_{n=1}\frac{q^{n^2}}{(q,q^2)_n}=\frac{q}{\eta\left(q^4\right)}\sum^{\infty}_{n=-\infty}\frac{(-1)^nq^{6n^2+6n}}{1-q^{4n+1}}.
\end{equation}
Then
$$
\psi(-q)=\frac{-q^{1/2}}{2\eta\left(q^4\right)}\sum^{\infty}_{n=-\infty}\frac{(-1)^n q^{6n^2+4n}}{\cosh\left(2\pi i z (2n+1/2)\right)}=
$$
$$
-\frac{q^{1/2}}{8\pi iz\cdot\eta\left(q^4\right)}\int^{+\infty}_{-\infty}\vartheta_4\left(\gamma,e\left(6 z\right)\right)e^{-i\gamma/2}\sec\left(\frac{\gamma}{4z}+\pi\right)d\gamma.
$$
Hence we can write the mock theta function as
\begin{equation}
\psi(-q)=\frac{q^{1/2}}{8\pi i z\cdot \eta\left(q^4\right)}\int^{+\infty}_{-\infty}\vartheta_4\left(\gamma,e(6z)\right)e^{-i\gamma/2}\sec\left(\frac{\gamma}{4z}\right)d\gamma.
\end{equation}
\\

The next theorem is paralel generalization of Theorem 9.\\
\\
\textbf{Theorem 18}\\
Let $q=e(z)$, $Im(z)>0$, $Im(w)\neq 0$. If we set 
$$
\widetilde{P}_j(a,b;A,B;z,w;\sigma):=
$$
\begin{equation}
=\frac{e^{-2iBb\pi z/A}}{2A\pi i w}\int^{+\infty \sigma}_{-\infty \sigma}\vartheta_3\left(\gamma,e(az)\right)e^{i\gamma^2 j/(2\pi az)}e^{-2iB\gamma/A}\sec\left(\frac{\gamma+b\pi z}{2Aw}\right)d\gamma,
\end{equation}
for $j=0,1$, then if $a>0$ we have
\begin{equation}
\widetilde{P}_0(a,b;A,B;z,w;\sigma)=\sum^{\infty}_{n=-\infty}\frac{q^{an^2+bn}}{\cosh(2\pi i w(An+B))}\textrm{, }\sigma\in\textbf{R}^{*}_{+}
\end{equation}
and holds
\begin{equation}
\widetilde{P}_j\left(a',b';A,B';z',w';z'\right)=\sqrt{-2iaz}\cdot\widetilde{P}_{1-j}\left(a,b;A,B;z,w;a\right),
\end{equation}
where the $a',b',B',z',w'$ are given from (143) below.\\
\\
\textbf{Notes.}\\
1) From (132) using (123) we can write
$$
\frac{1}{\eta_D(2az)}\sum^{\infty}_{n=-\infty}\frac{q^{an^2+bn}}{\cosh(2\pi i w(An+B))}=
$$
\begin{equation}
=\frac{ze^{-2iBb\pi z/A}}{2A i w}\int^{+\infty/z}_{-\infty/z}F_3\left(a,\gamma,z\right)e^{-2i\pi B\gamma z/A}e^{-i\pi \gamma^2 z/(2a)}\sec\left(\frac{(\gamma+b)\pi z}{2Aw}\right)d\gamma.
\end{equation}
2) From the deffinition relation (91) and from (130) of Theorem 17 we get
\begin{equation}
\sum^{\infty}_{n=-\infty}\frac{q^{an^2+bn}}{\cosh(2\pi i w (An+B))}=\frac{e^{-2iBb\pi z/A}}{A}P_0\left(a,b;z,wA;\frac{-2B}{A}\right).
\end{equation}
3) If we set
$$
P_G(a,b;A,B;z,w;x):=
$$
\begin{equation}
\frac{e^{-2iBb\pi z/A}}{2A\pi i w}\int^{+\infty}_{-\infty}\vartheta_3\left(\gamma,e(az)\right)e^{-2iB\gamma/A}\sec\left(\frac{\gamma+b\pi z}{2Aw}\right)e^{i\gamma x}d\gamma,
\end{equation}
then
\begin{equation}
P_G(a,b;A,B;z,w;0)=\sum^{\infty}_{n=-\infty}\frac{q^{an^2+bn}}{\cosh(2\pi i w(An+B))}.
\end{equation}
The $P_G(a,b;A,B;z,w;x)$ is a generalization of this Lerch sum (for $x=0$ it becomes the Lerch sum (141)). The Fourier transform of $P_G(a,b;A,B;z,w;x)$ with respect to $x$ is a known function:
\begin{equation}
P_G(a,b;A,B;z,w;(.))\wedge(s)
=\frac{e^{-2\pi i Bb z/A}}{iAw}\vartheta_3\left(s,e(az)\right)e^{-2iBs/A}\sec\left(\frac{s+b\pi z}{2Aw}\right).
\end{equation}
If we consider the transformation of variables
\begin{equation}
a'=\frac{1}{a}\textrm{, }b'=2bza'\textrm{, }z'=\frac{-1}{4z}\textrm{, }B'=-2B a z\textrm{, }w'=2wa'z'\textrm{, }s'=2sa'z',
\end{equation}
then from (84),(85) we get
$$
P_G(a',b';A,B;z',w';(.))\wedge(s')=
$$
$$
\frac{e^{-2\pi i Bb'z'/A}}{iAw'}\vartheta_3\left(s',e(a'z')\right)e^{-2iBs'/A}\sec\left(\frac{s'+b'\pi z'}{2Aw'}\right)=
$$
$$
-2iaz\frac{e^{\pi i Bb/(aA)}}{-Aw}\sqrt{-2iaz}\vartheta_3\left(s,e(az)\right)e^{is^2/(2\pi a z)}e^{iBs/(azA)}\sec\left(\frac{s}{2Aw}+\frac{bz\pi }{2Aw}\right).
$$
Hence
$$
P_G(a',b';A,B';z',w';(.))\wedge(s')=
$$
$$
-i(-2iaz)^{3/2}\frac{e^{-2\pi i Bbz/A}}{iAw}\vartheta_3\left(s,e(az)\right)e^{is^2/(2\pi a z)}e^{-2iBs/A}\sec\left(\frac{s}{2Aw}+\frac{bz\pi }{2Aw}\right).
$$
Hence we have the next\\
\\
\textbf{Theorem 19.}\\
If $a>0$, $Im(z)>0$, $Im(w)\neq 0$, $A>0$, and $a',b',z',w',s'$ as in (143), then
$$
P_G\left(a',b';A,B';z',w';(.)\right)\wedge\left(s'\right)=
$$
$$
=-i(-2iaz)^{3/2}\exp\left(\frac{is^2}{2\pi a z}\right)P_G(a,b;A,B;z,w;(.))\wedge(s).
$$
\\

We proceed by setting
\begin{equation}
S_0(a,b;A,B;z,w):=\frac{1}{\eta_D(2az)}\sum^{\infty}_{n=-\infty}\frac{q^{an^2+bn}}{\cosh(2\pi i w(An+B))}.
\end{equation}
Then from (132) we have with change of variable $\gamma\rightarrow \gamma z$, ($a>0$):
$$
S_0(a,b;A,B;z,w)=
$$
$$
\frac{e^{-2iBb\pi z/A}z}{2A i w}\int^{+\infty/z}_{-\infty/z}F_3\left(\frac{1}{a},\frac{2\gamma z}{a};\frac{-1}{4z}\right)e^{-2i\pi B\gamma z/A}\sec\left(\frac{(\gamma+b)\pi z}{2Aw}\right)d\gamma=
$$
$$
\frac{e^{-2iBb\pi z/A}z}{2A i w}\int^{+\infty/z}_{-\infty/z}F_3\left(a',\gamma';z'\right)e^{-2i\pi B\gamma z/A}\sec\left(\frac{(\gamma+b)\pi z}{2Aw}\right)d\gamma=
$$
$$
\frac{e^{-2iBb\pi z/A}z}{2A i w}\int^{+\infty/z}_{-\infty/z}F_3\left(a,\gamma;z\right)e^{-i\pi \gamma^2 z/(2a)}e^{-2i\pi B\gamma z/A}\sec\left(\frac{(\gamma+b)\pi z}{2Aw}\right)d\gamma.
$$
Hence
$$
S_0(a',b;A,B;z',w)=
$$
$$
\frac{e^{-2iBb\pi z'/A}z'}{2A i w}\int^{+\infty/z'}_{-\infty/z'}F_3\left(a',\gamma;z'\right)e^{-i\pi\gamma^2z'/(2a')}e^{-2i\pi B\gamma z'/A}\sec\left(\frac{(\gamma+b)\pi z'}{2Aw}\right)d\gamma=
$$
$$
\frac{e^{-2iBb\pi z'/A}z'2z}{2aA i w}\int^{+\infty/(z'z)}_{-\infty/(z'z)}F_3\left(a',2\gamma z/a;z'\right)e^{-i\pi4\gamma^2z^2z'/(2a^2a')}\times
$$
$$
\times e^{-2i\pi B2\gamma z z'/(aA)}\sec\left(\frac{(2\gamma z/a+b)\pi z'}{2Aw}\right)d\gamma=
$$
$$
\frac{-e^{iBb\pi/(2Az)}}{4aAiw}\int^{-\infty}_{+\infty}F_3(a',\gamma',z')e^{i\pi\gamma^2z/(2a)}e^{\pi i \gamma B/(aA)}\sec\left(\frac{\gamma\pi}{4Aaw}+\frac{b\pi}{8Azw}\right)d\gamma.
$$
Hence using (123), we get\\
\\
\textbf{Theorem 20.}\\
Let $a>0$, $Im(z)>0$, $Im(w)\neq 0$ and $A>0$, with $a',b',B',z',w'$ be as in (143). Then if we set
$$
\widetilde{S}_j(a,b;A,B;z,w;\sigma):=
$$
\begin{equation}
\frac{e^{-2\pi iBb z/A}}{2Aiw}\int^{+\infty\sigma}_{-\infty\sigma}F_3\left(\frac{1}{a},\frac{2\gamma}{a},\frac{-1}{4z}\right)e^{i\pi\gamma^2 j/(2 a z)}e^{-2\pi iB\gamma/A}\sec\left(\frac{\gamma\pi}{2Aw}+\frac{b\pi z}{2Aw}\right)d\gamma,
\end{equation}
it holds for $j=0,1$:
$$
\widetilde{S}_j(a',b';A,B';z',w';\sigma)=-\widetilde{S}_{1-j}(a,b;A,B;z,w;\sigma a z).
$$
Also if $\sigma>0$, then
\begin{equation}
\widetilde{S}_0(a,b;A,B;z,w;\sigma)=\frac{1}{\eta_D(2az)}\sum^{\infty}_{n=-\infty}\frac{q^{an^2+bn}}{\cosh(2\pi i w(An+B))}.
\end{equation}
\\

We can also generalize the function $S$ of (124) as folows:
$$
S_{G}(a,b;A,B;z,w;x):=
$$
\begin{equation}
\frac{e^{-2\pi i Bbz/A}}{2Aiw}\int^{+\infty}_{-\infty}F_3\left(\frac{1}{a},\frac{2\gamma}{a},\frac{-1}{4z}\right)e^{-2\pi i\gamma B/A}\sec\left(\frac{\gamma\pi}{2Aw}+\frac{b\pi z}{2Aw}\right)e^{i\gamma x}d\gamma.
\end{equation}
Hence
$$
S_G(a,b;A,B;z,w;0)=S_0(a,b;A,B;z,w)=
$$
\begin{equation}
=\frac{1}{\eta_D\left(2za\right)}\sum^{\infty}_{n=-\infty}\frac{q^{an^2+bn}}{\cosh\left(2\pi i w (An+B)\right)}.
\end{equation}
By this generalization we get the next formula evaluating the Fourier transform of $S_G$:
$$
S_{G}\left(a,b;A,B;z,w;(.)\right)\wedge(s)=
$$
\begin{equation}
=\frac{e^{-2\pi iBbz/A}}{2Aiw}F_3\left(\frac{1}{a},\frac{2s}{a},\frac{-1}{4z}\right)e^{-2\pi i s B/A}\sec\left(\frac{s\pi}{2Aw}+\frac{b\pi z}{2Aw}\right).
\end{equation}
\\
\textbf{Theorem 21.}\\
If $a>0$, $Im(z)>0$, $Im(w)\neq 0$ and $A>0$, then for the transformation of variables
\begin{equation}
a'=\frac{1}{a}\textrm{, }b'=2ba'z\textrm{, }B'=-2Baz\textrm{, }z'=\frac{-1}{4z}\textrm{, }w'=2wa'z'\textrm{, }s'=2sa'z',
\end{equation}
holds
\begin{equation}
S_G\left(a',b';A,B';z',w';(.)\right)\wedge(s')=-2az e^{\pi is^2/(2az)} S_G\left(a,b;A,B;z,w;(.)\right)\wedge(s).
\end{equation}
\\
\textbf{Proof.}
$$
S_G\left(a',b';A,B';z',w';(.)\right)\wedge(s')=
$$
$$
=\frac{e^{-2\pi iB'b'z'/A}}{2Aiw'}F_3\left(\frac{1}{a'},\frac{2s'}{a'},\frac{-1}{4z'}\right)e^{-2\pi i s' B'/A}\sec\left(\frac{s'\pi}{2Aw'}+\frac{b'\pi z'}{2Aw'}\right)=
$$
$$
=\frac{e^{-2\pi i (-2Bza)2ba'z z'/A}}{2Ai2wa'z'}F_3\left(a,\frac{4sa'z'}{a'},z\right)e^{-2\pi i 2sa'z'(-2Baz)/A}\times
$$
$$
\times\sec\left(\frac{2sa'z'\pi}{2A2wa'z'}+\frac{2ba'z\pi z'}{2A 2 w a'z'}\right)=
$$
$$
=-\frac{e^{-2\pi i Bbz/A}az}{Aiw}F_3\left(a,-\frac{s}{z},z\right)e^{-2\pi i sB/A}
\sec\left(\frac{s\pi}{2Aw}+\frac{bz\pi}{2A w }\right)=
$$
$$
=-\frac{e^{-2\pi i Bbz/A}az}{Aiw}F_3\left(\frac{1}{a},\frac{2s}{a},\frac{-1}{4z}\right)e^{i\pi s^2/(2az)}e^{-2\pi i sB/A}
\sec\left(\frac{s\pi}{2Aw}+\frac{bz\pi}{2A w}\right)=
$$
$$
=-2aze^{i\pi s^2/(2az)}S_G\left(a,b;A,B;z,w;(.)\right)\wedge(s).
$$
$qed$

\[
\]

\centerline{\bf References}

[1]: J.V. Armitage, W.F. Eberlein. ''Elliptic Functions''. Cambridge University Press. (2006)\\

[2]: N.D. Bagis. ''Some Results on Infinite Series and Divisor Sums''. (2009). Revised 2018. arXiv:0912.4815v3 [math.GM]\\ 

[3]: N.D. Bagis. ''On the Complete Evaluation of Jacobi Theta Functions''. (2015). Revised 2019. arXiv:1503.01141v2 [math.GM]\\

[4]: C.K. Chui. ''An Introduction to Wavelets''. Academic Press, Inc. 1992.\\

[5]: Don Zagier. ''Ramanujan's Mock Theta Functions and their Applications (after Zwegers and Ono-Bringmann)''. (2007), page stored in the Web.\\

[6]: Sander Zwegers. ''Mock Theta Functions''. Ph.D thesis (2002).\\arXiv:0807.4834v1\\ 

[7]: Byungchan Kim, Jeremy Lovejoy. ''Ramanujan-type partial theta identities and conjugate Bailey pairs, II. Multisums''. Ramanujan Journal. (2018). Vol 46. pages: 743-764.\\

[8]: D.R. Hickerson, E.T. Mortenson. ''Hecke-type double sums, Appell-Lerch sums, and mock theta functions (I)''. Proc. Lond. Math. Soc. (3) 109(2), 382-422 (2014).\\

[9]: E.T. Mortenson. ''On the dual nature of partial theta functions and Appell-Lerch sums''. Adv.Math. 264, 236-260 (2014).\\ 

[10]: K. Bringmann, K. Ono. ''The $f(q)$ mock theta function conjecture and partition ranks''. Inventiones Mathematicae. 165(2). pg 243-266,(2006).\\ 

[11]: K. Bringmann, K. Ono. ''Dyson's ranks and Maass forms''. Ann.Math. 171, pg 419-449. (2010).\\  

[12]: Eric M. Rains. ''Multivariate Quadratic Transformations and the Interpolation Kernel''. SIGMA 14 (2018), 019, 69 pages.\\

[13]: A. Dadholkar, S. Murthy, D. Zagier. ''Quantum Black Holes, Wall Crossing and Mock Modular Forms''. arXiv:1208.4074v2 [hep-th] 3 Apr 2014.

\end{document}